\documentclass[preprint,review,12pt]{elsarticle}

\usepackage{graphicx}

\usepackage{amssymb}
\usepackage{amsthm}

\usepackage{url}
\newtheorem{theorem}{Theorem}[section]
\newtheorem{lemma}{Lemma}[section]
\newtheorem{algorithm}{Algorithm}[section]
\newtheorem{remark}{Remark}[section]
\newtheorem{assumption}{Assumption}[section]

\journal{}

\begin{document}
\begin{frontmatter}


\title{A New Superlinearly Convergent Algorithm of Combining QP Subproblem with System of Linear
Equations for Nonlinear Optimization\tnoteref{s1}}
 \tnotetext[s1]{This work is supported by National
Natural Science Foundation of China (Grant Nos. 11071158, 11101262)}

\author[gxu]{Jin-Bao Jian\corref{cor2}}
\ead{jianjb@gxu.edu.cn}

\author[shu]{Chuan-Hao Guo\corref{cor1}}
\ead{guo-ch@live.cn, phone number:+86-1381-8217-112}

\author[gxu]{Chun-Ming Tang\corref{cor2}}
\ead{cmtang@gxu.edu.cn}

\author[shu]{Yan-Qin Bai}
\ead{yqbai@shu.edu.cn}

\cortext[cor1]{Corresponding author}

\address[gxu]{College of Mathematics and Information Science,
Guangxi University, Nanning, Guangxi 530004, China}

\address[shu]{Department of Mathematics, Shanghai University,
Shanghai 200444, China}

\begin{abstract}
In this paper, a class of optimization problems with nonlinear
inequality constraints is discussed. Based on the ideas of
sequential quadratic programming algorithm and the method of
strongly sub-feasible directions, a new superlinearly convergent
algorithm is proposed. The initial iteration point can be chosen
arbitrarily for the algorithm. At each iteration, the new algorithm
solves one quadratic programming subproblem which is always
feasible, and one or two systems of linear equations with a common
coefficient matrix. Moreover, the coefficient matrix is uniformly
nonsingular. After finite iterations, the iteration points can
always enter into the feasible set of the problem, and the search
direction is obtained by solving one quadratic programming
subproblem and only one system of linear equations. The new
algorithm possesses global and superlinear convergence under some
suitable assumptions without the strict complementarity. Finally,
some preliminary numerical experiments are reported to show that the
algorithm is promising.
\end{abstract}


\end{frontmatter}
\section{Introduction}\label{introduction}
In this paper, we consider the following nonlinear inequality
constrained optimization problem
$$\rm{(NCP)}\ \ \ \ \
\begin{array}{lll}
&\min\ \ f_0(x)\\
&{\rm s.t.} \ \ \ f_j(x)\leq 0,\ \ j\in I\triangleq\{1,2,\dots,m\},
 \end{array}
 $$
  where $x\in{R^n}$ and the functions
 $f_{j}(x): R^{n} \rightarrow R \ (j\in\{0\}\cup I)$ are all
 continuously differentiable. We denote the feasible set and gradients for problem (NCP) as
follows
$$X=\{x\in R^n:\ f_{j}(x)\leq 0, j\in I\},\ \ \ \ g_j(x)=\nabla f_j(x),\ j\in\{0\}\cup I.$$

It is well-known that sequential quadratic programming (SQP)
algorithms are acknowledged to one of successful algorithms
available for solving problem (NCP) and have good superlinear
convergence properties, they have been widely studied and
investigated by many authors in \cite{pt1987,bt1995,gt2000,gr2010}.

The iterative process of the standard SQP algorithms is as follows.
Let the current iteration point be $x$. Computing a search direction
$\bar{d}$ by solving the following quadratic programming (QP)
subproblem $$\ \ \ \
\begin{array}{lll}
&\min\ \ g_0(x)^Td+\frac{1}{2}d^TBd\\
&{\rm s.t.} \ \ \ f_j(x)+g_j(x)^Td\leq 0,\ \ j\in I,
 \end{array}
 $$
where $B\in R^{n\times n}$ is a symmetric matrix that approximates
the Hessian of the Lagrangian function associated with problem (NCP)
at $(x,\lambda)$, and $\lambda$ is a vector of nonnegative Lagrange
multiplier estimates. Perform a search to determine a steplength $t$
and let the next iteration point be $\bar{x}=x+t\bar{d}$.

However, in the classical SQP algorithms, there are two
shortcomings: (i) (QP) subproblem may be inconsistent, i.e., the
feasible set of (QP) subproblem may be empty; (ii) the Maratos
effect \cite{m1978} may occur, i.e., a full step of one can not be
taken close to a solution of problem (NCP). In order to overcome
disadvantage (i), various techniques have been proposed in
\cite{pt1987,ly2000}. A popular way to overcome disadvantage (ii) is
to use a higher-order direction, which is generated by solving a
(QP) subproblem \cite{pt1987} or a system of linear equations (SLE)
\cite{zz2001}, or directly given by an explicit formula
\cite{jzth2006}.

For SQP algorithms, feasible SQP (FSQP) algorithms are particularly
useful for solving those problems arising from engineering design
where the objective function $f_0$ might be undefined outside the
feasible set $X$. Another advantage of FSQP algorithms is that the
objective function $f_0$ can be used as a merit function to avoid
the use of a penalty function. In particular, Panier and Tits
\cite{pt1987} present a FSQP algorithm in which the generated
iteration points lie in the feasible set $X$. Two or three (QP)
subproblems need to be solved at each iteration. In order to obtain
the global convergence, they need to strengthen the requirement on
the first-order feasible descent condition. The superlinear
convergence rate is proved under the strict complementarity
assumption. Zhu and Jian \cite{zj2005} further improve the algorithm
in \cite{pt1987}. They introduce a new definition for the
first-order feasible condition which is weaker than the first-order
feasible descent condition in \cite{pt1987}, and propose a new FSQP
algorithm based on this new condition. The strict complementarity
assumption is also necessary for obtain the superlinear convergence.

One shortcoming of FSQP algorithms is that they are usually require
a feasible starting point, while computing such a point is generally
a nontrivial work \cite{bv2004}. In order to overcome this
shortcoming, Polak et al. \cite{ptm1979} propose a combined phase
I-phase II algorithm with arbitrary initial point for solving
problem (NCP). This algorithm becomes a method of feasible
directions (MFD) \cite{z1960} when iteration points enter into the
feasible set $X$. Jian further improve algorithm \cite{ptm1979} and
propose a method of strongly sub-feasible direction in \cite{j1995},
which not only unified automatically the operations of minimization
(Phase I), but also guaranteed that the number of the functions
satisfying the inequality constraints is nondecreasing. Since their
algorithms only using the information of first-order derivatives,
the algorithms in \cite{ptm1979,j1995} converge linearly at best.

In this paper, motivated by the ideas in \cite{zj2005,j1995}, we
propose a new algorithm combining (QP) subproblem with method of
strongly sub-feasible directions for solving problem (NCP). Unlike
algorithm in \cite{zj2005}, a descent direction $d_0$ for $f_0$ at
$x$ is obtained by solving a (QP) subproblem which is always
feasible. For obtaining the global convergence, a mere feasible
direction $\tilde{d}$ is obtained by solving a SLE. Then, $d_0$ is
tilted by making a convex combination
$q=(1-\beta)d_0+\beta\tilde{d}$ of $d_0$ and $\tilde{d}$, where
$\beta$ satisfies certain condition. With the help of this convex
combination, the global convergence is proved. In order to overcome
the Maratos effect and obtain the superlinear convergence, a
higher-order correction direction $d_1$ is computed by a new SLE
which has a common coefficient matrix with the previous SLE. Under
the strong second-order sufficient conditions without the strict
complementarity, the new algorithm is proved to be superlinearly
convergent. Moreover, the initial iteration point is chosen
arbitrarily, and after finite iterations, the iteration points can
always enter into the feasible set $X$. Finally, ,some numerical
results are reported to shown that the proposed algorithm is
promising.

At the end of this section, the main features of the proposed
algorithm are summarized as follows:

$\bullet$ the initial iteration point is arbitrary, and the number
of satisfied constraint functions is nondecreasing;

$\bullet$ the objective function of problem (NCP) is used directly
as the merit function, and the line search techniques are different
from others;

$\bullet$ at each iteration, the search direction is generated by
solving one QP subproblem and one or two SLEs with the same
coefficient matrix;

$\bullet$ after finite iterations, the search direction is obtained
by solving one QP subproblem and only one SLE, and the iteration
points always enter into the feasible set $X$;

$\bullet$ under some mild conditions without the strict
complementarity, the proposed algorithm possesses global and
superlinear convergence.

This paper is organized as follows. In the next Section
\ref{section2}, we present the details of our algorithm and discuss
its properties. In Sections \ref{section3} and \ref{section4}, the
algorithm is proved to possess global and superlinear convergence,
respectively. In Section \ref{section5}, some elementary numerical
experiments are reported. Finally, some concluding remarks are given
in Section \ref{section6}.
\section{Description of algorithm}\label{section2}
To simplify the analysis, we use the following notations
\begin{equation}\label{1}I^{-}(x)=\{j\in{I}: f_{j}(x)\leq 0\},\ \
I^{+}(x)=\{j\in{I}: f_{j}(x)>0\},\end{equation}
 \begin{equation}\label{2}
\varphi(x)=\max\{0,\ f_{j}(x),\ j\in{I}\}=\max\{0,\ f_{j}(x),\ j\in
I^{+}(x)\},\end{equation}
 \begin{equation}\label{3} \bar{f}_j(x)=\left\{\begin{array}{ll}f_j(x),
&j\in
I^-(x);\\
f_j(x)-\varphi(x), &j\in I^+(x),\end{array}\right. \ \ I(x)=\{j\in
I:\ \bar{f}_{j}(x)=0\}.\end{equation}

Assume that the following two basic assumptions for problem (NCP)
are hold throughout this paper.
\begin{assumption}\label{assumption1} The functions $f_j\ (j\in
\{0\}\cup I)$ are all first-order continuously differentiable.
\end{assumption}
\begin{assumption}\label{assumption2} The
gradient vectors $\{g_{j}(x): j\in I(x)\}$ are linearly independent
for each $x\in R^n$.
\end{assumption}
\begin{remark}\label{remark1}
In Assumption \ref{assumption2}, the linearly independent gradients
contains two parts: the one part is the gradients of the feasible
constraint functions in active set, and the other part is the
gradients of the maximal violated functions. This assumption is
becoming the standard linearly independent constraint qualification
(LICQ) only if the iteration point is feasible. Moreover, this
assumption plays a big role in the analysis of the following Lemmas
\ref{lemma2}, \ref{lemma3} and Theorem \ref{theorem1}.
\end{remark}

For the current iteration point $x^k$ and an associated symmetric
positive definite matrix $B_k$, using the notations above, we
introduce the following (QP) subproblem \cite{j2000}
$$\rm{(QPs)}\ \ \ \ \
\begin{array}{lll} &\min\ \ g_0(x^k)^Td+\frac{1}{2}d^TB_kd
\\ &{\rm s.t.} \ \ \
\bar{f}_j(x^k)+g_j(x^k)^Td\leq 0,\ \ \forall j\in I,\\
  \end{array}$$
and we denote simply
$$I_k^-=I^-(x^k),\ I_k^+=I^+(x^k),\ I_{k}=I(x^k),\ \varphi_k=\varphi(x^k).$$
It is obviously that subproblem (QPs) has the following merits:\\
(i) subproblem (QPs) is always feasible with feasible solution
$d=0$;\\
(ii) subproblem (QPs) is a strictly convex program while $B_k$ is
positive definite, so it has an (unique) optimal solution;\\
(iii) $d$ is a solution of subproblem (QPs) if and only if it is a
KKT point of subproblem (QPs).

Let $d_0^k$ is an optimal solution of subproblem (QPs) at the $k-th$
iteration, i.e., there exists a
 corresponding KKT multiplier vector $\lambda^k=(\lambda_j^k,\ j\in I)$
 such that
 \begin{equation}\label{4} \left\{\begin{array}{ll} g_0(x^k)+B_kd_0^k+\sum\limits_{j\in
 I}\lambda_j^kg_j(x^k)=0,\\
  \bar{f}_j(x^k)+g_j(x^k)^Td_0^k\leq 0,\ \lambda_j^k\geq 0,\
 \lambda_j^k(\bar{f}_j(x^k)+g_j(x^k)^Td_0^k)=0,\ \forall j\in
 I.
 \end{array}\right.
 \end{equation}
From (\ref{4}) and the KKT condition for problem (NCP), the
following lemma holds immediately.
\begin{lemma}\label{lemma1} $x^k$
 is a KKT point for problem (NCP) if and only if $(d_0^k,\ \varphi_k)=(0,\ 0)$.
\end{lemma}

On one hand, in view of $d=0$ is a feasible solution and $d_0^k$ is
an optimal solution for subproblem (QPs), respectively, it follows
that
\begin{equation}\label{5}g_0(x^k)^Td_0^k+\frac{1}{2}(d_0^k)^TB_kd_0^k\leq
0,\end{equation} and if $d_0^k\neq 0$, together with the positive
definiteness of $B_k$, (\ref{5}) imply that
$$g_0(x^k)^Td_0^k\leq -\frac{1}{2}(d_0^k)^TB_kd_0^k<0,$$ i.e., $d_0^k$
is a descent direction for problem (NCP) at $x^k$.

On the other hand, $d_0^k$ may be not a feasible direction for
problem (NCP) at the feasible iteration point $x^k$. Even when
$d_0^k$ is a feasible direction, a line search may not allow a full
step of one to be taken in a neighborhood of an optimal solution and
thus superlinear convergence may never take place. In order to get
an improving direction and taking into account that $x^k$ may be
infeasible, we first propose a new SLE
\begin{equation}\label{6} V_k\left(\begin{array}{c}d
\\ h\end{array}\right)\triangleq\left(\begin{array}{c}B_k\ \ \ N_k\\N_k^T\ \
-D^k\end{array}\right)\left(\begin{array}{c}d\\h
\end{array}\right)=\left(\begin{array}{c}0\\-(||d_0^k||+\varphi_k^\sigma)\varpi\end{array}\right)
\end{equation}
 to generate an updated direction $\tilde{d}^k$, where $0\in R^n$,
$\varpi=(1,1,\ldots,1)^T\in R^m$, $\sigma\in(0,1)$ and
\begin{equation}\label{7} \left\{\begin{array}{ll}N_k=(g_{j}(x^k),\ j\in I),\ \ \ \
D^k={\rm diag}(D_{j}^k,\ j\in I),\\
D_j^k=|\bar{f}_j(x^k)|(|\bar{f}_j(x^k)+g_j(x^k)^Td_0^k|+||d_0^k||),\
j\in I.\end{array}\right.\end{equation} From (\ref{6}) and
(\ref{7}), it follows that
$$g_j(x^k)^T\tilde{d}^k=-\|d_0^k\|-\varphi_k^\sigma,\ \forall j\in I_k,$$
i.e., $\tilde{d}^k$ is a mere feasible direction.

The following lemma describes the solvability of SLE (\ref{6}), its
proof is similar to Lemma 2.2 in \cite{jkzt2009} and is omitted
here.
\begin{lemma}\label{lemma2} Suppose that Assumptions \ref{assumption1} and \ref{assumption2}
 hold and $B_k$ is positive definite. Then, the coefficient
matrix $V_k$ defined in (\ref{6}) is nonsingular and (\ref{6}) has a
unique solution. \end{lemma}

Then, in order to yield an improving search direction at iteration
point $x^k$ (feasible or infeasible), we consider a convex
combination of $d_0^k$ and $\tilde{d}^k$ as follows
\begin{equation}\label{8}
q^k=(1-\beta_k)d_0^k+\beta_k\tilde{d}^k,\end{equation} where
$\beta_k$ is the maximum value of $\beta\in[0,1]$ that satisfies the
following relationship
\begin{equation}\label{9}g_0(x^k)^Tq^k=(1-\beta_k)g_0(x^k)^Td_0^k+\beta_kg_0(x^k)^T\tilde{d}^k\leq
\theta g_0(x^k)^Td_0^k+\varphi_k^\theta,\end{equation} where the
positive parameter $\theta<\sigma$. Taking into account that
$g_0(x^k)^Td_0^k\leq 0$, it follows that $\beta_k$ can be yielded by
the following explicit formula
\begin{equation}\label{10}\beta_k=\left\{\begin{array}{ll}\min\{1,\
\frac{(\theta-1)g_0(x^k)^Td_0^k+\varphi_k^\theta}{g_0(x^k)^T\tilde{d}^k-g_0(x^k)^Td_0^k}\},&{\rm{if}}\
g_0(x^k)^T\tilde{d}^k>g_0(x^k)^Td_0^k;\\ 1,&{\rm{if}}\
g_0(x^k)^T\tilde{d}^k\leq
g_0(x^k)^Td_0^k.\end{array}\right.\end{equation}

\begin{remark}\label{remark2}From (\ref{9}), it follows that the rate of increase of
the objective function $f_0$ at point $x_k$ along direction $q^k$ is
bounded from upper by $\theta g_0(x^k)^Td_0^k+\varphi_k^\theta$, and
it is just so ensuring the direction $q^k$ possesses excellent
convergence. In addition, the parameter $\theta\in(0,\sigma)$ plays
important roles in avoiding Maratos effect for  analyzing the
request (\ref{14}) in the algorithm  as well as forcing the
iteration points always get into the feasible region after a finite
number of iterations. These can be seen in the latter analysis,
e.g., Lemma \ref{lemma6} and Theorems \ref{theorem3}, \ref{theorem4}
as well as \ref{theorem5}.
\end{remark}

From the relationships of (\ref{8}), (\ref{9}) and (\ref{6}), the
following lemma can be proved easily.
\begin{lemma}\label{lemma3} Suppose that
assumptions in Lemma \ref{lemma2} hold. Then\\ (i)
$g_0(x^k)^Tq^k\leq-\frac{1}{2}\theta(d_0^k)^TB_kd_0^k+\varphi_k^\theta$;\
(ii) $g_j(x^k)^Tq^k\leq -\beta_k(||d_0^k||+\varphi_k^\sigma),\
\forall j\in I_k.$ \end{lemma}

From Lemma \ref{lemma3}, it holds that $q^k$ is an improving
direction either for problem (NCP) or for the maximal violated
constrained function $\varphi(x)$. In order to overcome the
possibility of Maratos effect, a suitable higher-order correction
direction must be introduced by an appropriate approach.
Additionally, taking into account avoiding the strict
complementarity condition, we introduce another SLE
 \begin{equation}\label{11} V_k\left(\begin{array}{c}d\\h
\end{array}\right)=\left(\begin{array}{c}0\\-(||d_0^k||^\tau+\varphi_k^\sigma)\varpi-\tilde{F}
(x^k+d_0^k)\end{array}\right) \end{equation} to yield a higher-order
correction direction $d_1^k$, where $\tau\in(2,3)$ and
\begin{equation}\label{12} \ \
 \tilde{F}(x^k+d_0^k)=(f_j(x^k+d_0^k)-f_j(x^k)-g_j(x^k)^Td_0^k,\ j\in
 I).\end{equation}

\begin{remark}\label{remark3} The right-hand-side of (\ref{11}) (in particular, the
introducing of $\tilde{F}(x^k+d^k_0)$) as well as $d_1^k$ play an
important role in the discussion of superlinear convergence in
Section \ref{section4} without the strict complementarity, these can
be found in the proofs of Lemma \ref{lemma8} and Theorem
\ref{theorem3}. In traditional analysis \cite{pt1987,zj2005,lt2000},
$d_1^k$ must satisfy $||d_1^k||=O(||d_0^k||^2)$ that is called
second-order correction. In this paper, the term $\varphi_k$ is
added, so the relation between $d_1^k$ and $d_0^k$ will be different
from the traditional forms, see Lemma \ref{lemma8}(i).
\end{remark}

Now, based on the analysis above, we can present our algorithm as
follows.
\begin{algorithm}\label{algorithm1}\end{algorithm}

Parameters: $\gamma,\eta,\varepsilon\in(0,1),\ 0<\theta<\sigma<1,\
0<\varrho<\sigma,\ \xi,\zeta>0,\ \alpha\in(0,0.5),\ \rho>1,\
\delta>2,\ \tau\in(2,3).$

Data: $x^0\in R^n$, a symmetric positive definite matrix $B_0\in
R^{n\times n}.$ Set $k:=0$.

\textbf{Step 1.} Solve subproblem (QPs) to get a (unique) solution
$d_0^k$ and the corresponding KKT multiplier vector
$\lambda^k=(\lambda_j^k,\ j\in I)$. If $(d_0^k,\ \varphi_k)=0$, then
$x^k$ is a KKT point for problem (NCP) and stop.

\textbf{Step 2.} Compute the correction direction $d_1^k$ by solving
SLE (\ref{11}) with a solution $(d_1^k, h_1^k)$, and let
$d^k=d_0^k+d_1^k$. If
\begin{equation}\label{13}g_0(x^k)^Td_0^k\leq\zeta\min\{-||d_0^k||^\delta,\
-||d^k||^\delta\}+\xi\varphi_k^\varrho,\end{equation}
 then go to Step 3;
otherwise, go to Step 4.

\textbf{Step 3.} Let $t=1$,

\rm(a) if \begin{equation}\label{14}\left\{\begin{array}{ll}
f_0(x^k+td^k)\leq
f_0(x^k)+\alpha t g_0(x^k)^Td_0^k+\rho(1-\alpha)t\varphi_k^\theta,\\
  f_j(x^k+td^k)\leq \varphi_k-\alpha t(||d_0^k||^\tau+\varphi_k^\sigma),\ \ j\in
 I_k^+,\\
  f_j(x^k+td^k)\leq 0,\ \ j\in
 I_k^-,\\
 \end{array}\right.
 \end{equation}
 is satisfied for the current value $t$, then let $t_k=t$, go to Step 6; otherwise, go to (b) below.

\rm(b) Reset $t:=\frac{1}{2}t.$ If $t<\varepsilon$, then go to Step
4;
 otherwise, repeat part (a).

\textbf{Step 4.} Solve SLE (\ref{7}) to get $(\tilde{d}^k,
\tilde{h}^k)$, and compute $\beta_k$ according to (\ref{10}), then
obtain the direction $q^k$ by (\ref{8}).

\textbf{Step 5.} Compute the steplength $t_k$ which is the first
number $t$ in the sequence $\{1,\eta,\eta^2,\dots\}$ that satisfies
the following inequalities
\begin{equation}\label{15} f_0(x^k+tq^k)\leq f_0(x^k)+\gamma t
g_0(x^k)^Tq^k+\rho(1-\gamma)t\varphi_k^\theta,\end{equation}
\begin{equation}\label{16}f_j(x^k+tq^k)\leq \varphi_k-\gamma
t\beta_k(||d_0^k||+\varphi_k^\sigma),\ \ j\in
 I_k^+,\end{equation}
\begin{equation}\label{17}f_j(x^k+tq^k)\leq 0,\ \ j\in
 I_k^-,\end{equation}
 then let $d^k=q^k$.

\textbf{Step 6.} Compute a new symmetric positive definite matrix
$B_{k+1}$ by some suitable techniques, set $x^{k+1}=x^k+t_kd^k,\
k:=k+1,$ and go back to Step 1.

\begin{remark}\label{remark4} From the mechanism of
Algorithm \ref{algorithm1}, we know that when $\varphi_k\neq 0$, the
value of the maximal violated constrained functions (namely
$\varphi_k$) is strictly decreasing, moreover, the amounts of
decreasing have a close relation with $||d_0^k||+\varphi_k$. Thus,
the iteration points can approach the feasible region as soon as
possible. In particular, they can always getting into the feasible
set after a finite number of iterations.
\end{remark}

\begin{remark}\label{remark5} It is known that the role of the request (\ref{13}) is to
restrict the increasing  rate of the objective function $f_0$ at
point $x^k$ along direction $d_0^k$. But, theoretically speaking,
under the uniformly positive definite assumption (see Assumption
\ref{assumption3} below) on the sequence $\{B_k\}$ of matrices, the
request (\ref{13}) does not influence any theory analysis　of
Algorithm  \ref{algorithm1}, this can be seen in the latter
analysis. However, the request (\ref{13}) still has some influence
on the numerical effect of Algorithm  \ref{algorithm1}. From the
process of numerical experiments, it seems to be that, for
small-scale problems, the numerical results of Algorithm
\ref{algorithm1} with the request (\ref{13}) are better than the
case of ignoring this request, and the case is inverse for
middle-large-scale problems. In addition, the role of the exponents
$\delta>2$ and $\varrho<\sigma$ is to ensure, under suitable
assumptions, that the request (\ref{13}) is always satisfied when
the iteration point $x^k$ is close sufficiently to a KKT point, see
Lemma \ref{lemma8}(iii).
\end{remark}

\begin{remark}\label{remark6} There are two cycles between Step 1 and Step 6, i.e., cycle I:
Steps 1-2-3-6, and cycle II: Steps 1-2-4-5-6. Obviously, if the
algorithm successfully perform cycle I, the cost of computation is
relatively small. Fortunately, under suitable assumptions, we can
prove that Algorithm \ref{algorithm1} always performing cycle I
after a finite number of iterations.
\end{remark}

The lemma given below indicates that the line search in Step 5 of
Algorithm \ref{algorithm1} is well defined.

\begin{lemma}\label{lemma4} Suppose that Assumptions \ref{assumption1} and \ref{assumption2} hold.
Then, if Algorithm \ref{algorithm1} does not stop at Step 1, i.e.,
$(d_0^k,\ \varphi_k)\neq(0,\ 0)$, the line search
(\ref{15})-(\ref{17}) can be terminated after a finite number of
computations.
\end{lemma}
\begin{proof} Suppose that $(d_0^k,\ \varphi_k)\neq(0,\ 0)$ at the $k$-$th$
iteration. From (\ref{10}), it follows that $\beta_k>0$.

(1) Analyze the inequality (\ref{15}). From Taylor expansion and
Lemma \ref{lemma3}(i), it follows that
\begin{equation}\label{18}\begin{array}{ll}a_k(t)&\triangleq f_0(x^k+tq^k)-f_0(x^k)
-\gamma tg_0(x^k)^Tq^k-\rho(1-\gamma)t\varphi_k^\theta\\
\ \ \ &=(1-\gamma)tg_0(x^k)^Tq^k-\rho(1-\gamma)t\varphi_k^\theta+o(t||q^k||)\\
\ \ \
&\leq-\frac{1}{2}\theta(1-\gamma)t(d_0^k)^TB_kd_0^k-(\rho-1)(1-\gamma)t\varphi_k^\theta
+o(t||q^k||).\end{array}\end{equation} This together with
$\theta,\gamma\in(0,1),\ \rho>1$ and
$(d_0^k)^TB_kd_0^k+\varphi_k^\theta>0$ shows that $a_k(t)\leq 0$
holds for $t>0$ sufficiently small.

(2) Analyze the inequalities (\ref{16}).

(i) For $j\in I_{k}^+\cap I_k$, it follows that
$f_j(x^k)=\varphi_k$, expanding $f_j(x^k+tq^k)$ around $x^k$, and
combining Lemma \ref{lemma3}(ii), for $\gamma\in (0,1)$ and $t>0$
sufficiently small, we obtain
$$\begin{array}{ll}f_j(x^k+tq^k)&-\varphi_k+\gamma
t\beta_k(||d_0^k||+\varphi_k^\sigma)\\
&=tg_j(x^k)^Tq^k+\gamma t\beta_k(||d_0^k||+\varphi_k^\sigma)+o(t||q^k||)\\
&\leq-(1-\gamma)t\beta_k(||d_0^k||+\varphi_k^\sigma)+o(t||q^k||)\\
&\leq 0.\end{array}$$

(ii) For $j\in I^+_k\backslash I_k$, we have $f_j(x^k)-\varphi_k<0$,
and furthermore,
$$\lim\limits_{t\rightarrow 0^+}(f_j(x^k+tq^k)-\varphi_k+\gamma t\beta_k(||d_0^k||+\varphi_k^\sigma))
=f_j(x^k)-\varphi_k<0,$$ which implies the inequalities (\ref{16})
hold for $t>0$ sufficiently small.

(3) Analyze the inequalities (\ref{17}).

(i) For $j\in I_{k}^-\cap I_k$, expanding $f_j(x^k+tq^k)$ at $x^k$,
and taking into account Lemma \ref{lemma3}(ii), for $t>0$
sufficiently small, it follows that
$$f_j(x^k+tq^k)=tg_j(x^k)^Tq^k+o(t||q^k||)\leq -t\beta_k(||d_0^k||+\varphi_k^\sigma)+o(t||q^k||)\leq 0.$$

(ii) For $j\in I^-_k\backslash I_k$, we have $f_j(x^k)<0$. So, from
Assumption \ref{assumption1}, it follows that $f_j(x^k+tq^k)\leq 0$,
for $t>0$ sufficiently small.

Summarizing the analysis above, we conclude that there exits a
$\bar{t}_k>0$ such that the line search (\ref{15})-(\ref{17})
satisfies for all $t\in(0,\ \bar{t}_k]$ and the given conclusion
holds.
\end{proof}

At the end of this section, based on the line search conditions
(\ref{14}), (\ref{16}) and (\ref{17}), we can easily get the
following lemma.
\begin{lemma}\label{lemma5} Suppose that Assumptions
\ref{assumption1} and \ref{assumption2} hold. Then,\\ (i) for each
$k$, $I_k^-\subseteq I_{k+1}^-$ holds, so, if there exists an
iteration index $k_0$ such that $x^{k_0}\in X$, i.e.,
$\varphi_{k_0}=0$, then $x^k\in X$ for all $k\geq k_0$, and
$\{f_0(x^k)\}_{k\geq k_0}$ is decreasing;\\ (ii) if $x^k\not\in X$
and $x^{k+1}\not\in X$, then
$\varphi_{k+1}\leq\varphi_k-t_k\max\{\alpha(||d_0^k||^\tau+\varphi_k^\sigma),\
\gamma\beta_k(||d_0^k||+\varphi_k^\sigma)\}<\varphi_k;$ \\ (iii) for
$k$ large enough, the subsets $I_k^-$ and $I_k^+$ can be fixed,
i.e., $I_k^-\equiv I^-$ and $I_k^+\equiv I^+.$
\end{lemma}
\begin{remark}\label{remark7} From Lemma \ref{lemma5}(i) and (ii), it
follows that exactly one of the following two cases takes place:

Case A: There exists an iteration index $k_0$ such that
$\varphi_{k_0}=0$, then $\varphi_k\equiv0$ for all $k\geq k_0$;

Case B: For any $k=0,\ 1,\ 2,\ldots$, $\varphi_k>0$ and
$\varphi_{k+1}<\varphi_k$.
\end{remark}
\section{Global convergence}\label{section3}
In this section, we will establish the global convergence of
Algorithm \ref{algorithm1}. When Algorithm \ref{algorithm1} stops at
$x^k$, it follows that the iteration point $x^k$ is a KKT point for
problem (NCP) from Lemma \ref{lemma1} and Step 1 of Algorithm
\ref{algorithm1}. Now, suppose that an infinite sequence $\{x^k\}$
of iteration points is generated by Algorithm \ref{algorithm1}, and
we will show that every accumulation point $x^*$ of $\{x^k\}$ is the
KKT point for problem (NCP). For this purpose, the following
assumption is necessary.

\begin{assumption}\label{assumption3}The sequence $\{B_k\}$ of matrices is
uniformly positive definite, i.e., there exist two positive
constants $a$ and $b$ such that
\begin{equation}\label{19}a||d||^2\leq d^TB_kd\leq b||d||^2,\ \
\forall d\in R^n,\ \forall k.\end{equation}
\end{assumption}

Denote the active set for subproblem (QPs) by
\begin{equation}\label{20}J_k=\{j\in I:\
\bar{f}_j(x^k)+g_j(x^k)^Td_0^k=0 \}.\end{equation}

Suppose that $x^*$ is a given accumulation point of $\{x^k\}$. In
view of $I_k^+,\ I_k^-,\ J_k$ all being subsets of the finite set
$I$ and Lemma \ref{lemma5}(iii), we can assume without loss of
generality that there exists an infinite index set $K$ such that
\begin{equation}\label{21}x^k\rightarrow x^*,\ I_k^-\equiv I^-,\ I_k^+\equiv
I^+,\ J_k\equiv J,\ \varphi_k\rightarrow\varphi_*=\varphi(x^*),\
\forall k\in K.\end{equation} Based on the above conditions, we have
the following lemmas.
\begin{lemma}\label{lemma6} Suppose that Assumptions \ref{assumption1}, \ref{assumption2}
and \ref{assumption3} hold. Then\\
(i) sequence $\{d_0^k\}_{K}\triangleq\{d_0^k:\ k\in K\}$ is bounded;\\
(ii) there exists a constant $c>0$ such that $||V_k^{-1}||\leq c$
for all $k\in
K$, where $V_k$ is the coefficient matrix defined by (\ref{6}); \\
(iii) sequences $\{d_1^k\}_K$, $\{\tilde{d}^k\}_K$, $\{q^k\}_K$ and
$\{\tilde{h}^k\}_K$ are all bounded.
\end{lemma}
\begin{proof} (i) First, from $\lim\limits_{k\in K}x^k=x^*$ and the
continuity of $g_0(x)$, there exists a constant $\bar{c}>0$ such
that $||g_0(x^k)||\leq\bar{c}$ holds for all $k\in K$. Then, from
(\ref{5}) and Assumption \ref{assumption3}, it follows that
$-\bar{c}||d_0^k||+\frac{1}{2}a||d_0^k||^2\leq 0,\ \forall k\in K$,
which implies the boundedness of $\{d_0^k\}_K$.

(ii) The proof of conclusion (ii) is similar to Lemma 3.1 in
\cite{jkzt2009}, thus it is omitted here.

(iii) According to (\ref{6}), (\ref{11}) and parts (i) and (ii), it
follows that $\{d_1^k\}_K$, $\{\tilde{d}^k\}_K$ and
$\{\tilde{h}^k\}$ are all bounded. Furthermore, from (\ref{8}) and
$\beta_k\in[0,1]$, we obtain the boundedness of $\{q^k\}_K$.
\end{proof}

\begin{lemma}\label{lemma7} Suppose that assumptions stated in Lemma \ref{lemma6} hold.
Then $\lim\limits_{k\in K}(d_0^k,\varphi_k)=(0,0)$ and
$\lim\limits_{k\in K}\tilde{d}^k=\lim\limits_{k\in
K}q^k=\lim\limits_{k\in K}d_1^k=0$, where the index set $K$ is
defined by (\ref{21}).
\end{lemma}
\begin{proof} We first show that $\lim\limits_{k\in
K}(d_0^k,\varphi_k)=(0,0)$. Suppose by contradiction that there
exists an infinite index set $K^{'}\subseteq K$ and a constant
$\tilde{\epsilon}>0$ such that
$||(d_0^k,\varphi_k)||\geq\tilde{\epsilon},\ k\in K^{'}.$ So, there
exists an infinite set $K^{''}\subseteq K^{'}$ such that
$d_0^k\rightarrow d_0^*,\ \varphi_k\rightarrow\varphi_*,\
B_k\rightarrow B_*,\ k\in K^{''}$. Thus, in view of (\ref{5}),
(\ref{19}) and $\theta\in(0,1)$, we have
$$(\theta-1)g_0(x^k)^Td_0^k+\varphi_k^\theta\geq\frac{1}{2}(1-\theta)a||d_0^k||^2+\varphi_k^\theta
\stackrel{k\in
K^{''}}{\longrightarrow}\frac{1}{2}(1-\theta)a||d_0^*||^2+\varphi_*^\theta>0,$$
which further shows that
$$||d_0^k||^2+\varphi_k^\theta\stackrel{k\in
K{''}}{\longrightarrow}||d_0^*||^2+\varphi_*^\theta>0,\ \ \
||d_0^k||+\varphi_k^\sigma\stackrel{k\in
K^{''}}{\longrightarrow}||d_0^*||+\varphi_*^\sigma>0.$$ This
together with $(d_0^k, \varphi_k)\neq 0$ and (\ref{10}) implies that
there exists a constant $\bar{\epsilon}>0$ such that
\begin{equation}\label{22}\beta_k\geq\bar{\epsilon},\ \
||d_0^k||^2+\varphi_k^\theta\geq\bar{\epsilon}, \ \
||d_0^k||+\varphi_k^\sigma\geq\bar{\epsilon},\ \forall k\in
K^{''}.\end{equation} The rest proof is divided into two steps as
follows.

{\bf Step A}. Show that there exists a constant $\bar{t}>0$ such
that the steplength $t_k\geq\bar{t}$ for all $k\in K^{''}$.

From the mechanism of Algorithm \ref{algorithm1} we know that it is
sufficient to prove that the inequalities (\ref{15})-(\ref{17}) are
satisfied for all $k\in K^{''}$ and $t>0$ small enough.

(1) Analysis the inequality (\ref{15}). From the boundedness of
$\{d_0^k\}_K$ and (\ref{18}), (\ref{19}) as well as (\ref{22}), we
have
\begin{equation}\label{23}\begin{array}{ll}a_k(t)&\leq-(1-\gamma)t(\frac{1}{2}\theta(d_0^k)^TB_kd_0^k
+(\rho-1)\varphi_k^\theta)+o(t)\\
&\leq-(1-\gamma)t(\frac{1}{2}a\theta||d_0^k||^2+(\rho-1)\varphi_k^\theta)+o(t)\\
&\leq-(1-\gamma)t\min\{\frac{1}{2}a\theta,\
\rho-1\}(||d_0^k||^2+\varphi_k^\theta)+o(t)\\
&\leq-(1-\gamma)t\bar{\epsilon}\min\{\frac{1}{2}a\theta,\
\rho-1\}+o(t)\leq 0\end{array}\end{equation} holds for $t>0$ small
enough and all $k\in K^{''}$.

(2) Analysis the inequalities (\ref{16}).

(2-i) For $j\in I^+$ and $f_j(x^*)=\varphi_*$, from (\ref{7}) and
(\ref{3}), it follows that
$$D_j^k=|f_j(x^k)-\varphi_k|(|f_j(x^k)+g_j(x^k)^Td_0^k-\varphi_k|+||d_0^k||)\rightarrow0,\ k\in K^{''}.$$
Therefore, from Taylor expansion, (\ref{8}), (\ref{6}) and the
constraints of subproblem (QPs), for $t>0$ sufficiently small and
all $k\in K^{''}$, we obtain
$$\begin{array}{ll}&f_j(x^k+tq^k)-\varphi_k+\gamma
t\beta_k(||d_0^k||+\varphi_k^\sigma)\\
&=f_j(x^k)+tg_j(x^k)^Tq^k-\varphi_k+\gamma t\beta_k(||d_0^k||+\varphi_k^\sigma)+o(t)\\
&=f_j(x^k)-\varphi_k+t(1-\beta_k)g_j(x^k)^Td_0^k+t\beta_kg_j(x^k)^T\tilde{d}^k+\gamma
t\beta_k(||d_0^k||+\varphi_k^\sigma)+o(t)\\
&\leq(1-t(1-\beta_k))(f_j(x^k)-\varphi_k)-t\beta_k(||d_0^k||+\varphi_k^\sigma)
+t\gamma\beta_k(||d_0^k||+\varphi_k^\sigma)+o(t)\\
&\leq-(1-\gamma)t\beta_k(||d_0^k||+\varphi_k^\sigma)+o(t)\\
&\leq-\bar{\epsilon}^2(1-\gamma)t+o(t)\leq 0.\end{array}$$

(2-ii) For $j\in I^+$ and $f_j(x^*)<\varphi_*$, we have
$$\begin{array}{ll}f_j(x^k+tq^k)-\varphi_k+\gamma t\beta_k(||d_0^k||+\varphi_k^\sigma)&=f_j(x^k)-\varphi_k
+O(t)\\&\leq\frac{1}{2}(f_j(x^*)-\varphi_*)+O(t)\\&\leq
0\end{array}$$ holds for all $k\in K^{''}$ and $t>0$ sufficiently
small.

(3) Analysis the inequalities (\ref{17}). Similar to the analysis
for the inequalities (\ref{16}), one can show that (\ref{17}) hold
for all $k\in K^{''}$ and $t>0$ small enough.

Summarizing the analysis above, we conclude that there exists a
$\bar{t}>0$ such that $t_k\geq \bar{t}$ for all $k\in K^{''}$.

{\bf Step B.} Use $t_k\geq\bar{t}>0\ (k\in K^{''})$ to bring a
contradiction, and the discussion is divided into two cases.

Case I. Suppose that there exists an iteration index $k_0$ such that
$\varphi_{k0}=0$. Then $\{f_0(x^k)\}_{k\geq k_0}$ is decreasing.
Combining $\lim\limits_{k\in K^{''}}f_0(x^k)=f_0(x^*)$, it follows
that $\lim\limits_{k\rightarrow+\infty}f_0(x^k)=f_0(x^*)$. On the
other hand, taking into account the first inequality of (\ref{14}),
(\ref{15}), Lemma \ref{lemma3}(i), (\ref{19}) and (\ref{22}) as well
as $\varphi_k\equiv 0\ (\forall k\geq k_0)$, it follows that
$$\begin{array}{ll}f_0(x^{k+1})&\leq f_0(x^k)+t_k\max\{\alpha g_0(x^k)^Td_0^k,\ \gamma g_0(x^k)^Tq^k\}\\
&\leq f_0(x^k)+t_k\max\{-\frac{1}{2}\alpha a||d_0^k||^2,\
-\frac{1}{2}\theta\gamma a||d_0^k||^2\}\\
&\leq f_0(x^k)-\bar{t}\bar{\epsilon}\min\{\frac{1}{2}\alpha a,\
\frac{1}{2}\theta\gamma a\},\ \ \forall k\geq k_0,\ k\in
K^{''}.\end{array}$$
 Thus, passing to the limit $k\in K^{''}$ and
$k\rightarrow+\infty$ in the inequality above, we can bring a
contradiction.

Case II. Suppose that $\varphi_k>0$ for each $k$. Then
$\{\varphi_k\}_{k\geq 0}$ is decreasing. Combining
$\lim\limits_{k\in K^{''}}\varphi_k=\varphi_*$, one knows that
$\lim\limits_{k\rightarrow+\infty}\varphi_k=\varphi_*$. On the other
hand, taking into account the second inequality of (\ref{14}),
(\ref{16}), Lemma \ref{lemma3}(ii), (\ref{19}) and (\ref{22}) as
well as $\varphi_k>0\ (\forall k\in K^{''})$, it follows that
$$\begin{array}{ll}\varphi_{k+1}&\leq\varphi_k-t_k\max\{\alpha(||d_0^k||^\tau+\varphi_k^\sigma),
\ \gamma\beta_k(||d_0^k||+\varphi_k^\sigma)\}\\
&\leq\varphi_k-\bar{t}\bar{\epsilon}\max\{\alpha,\
\gamma\bar{\epsilon}\},\ \forall k\in K^{''}, k\ {\rm{large\
enough}}.\end{array}$$ Then, passing to the limit $k\in K^{''}$ and
$k\rightarrow+\infty$ in the above inequality, we can also bring a
contradiction.

Up to now, we have finished the proof of $\lim\limits_{k\in
K}(d_0^k,\varphi_k)=(0,0)$.

Finally, by Lemma \ref{lemma6}(ii), (\ref{6}), (\ref{11}), (\ref{8})
and $\lim\limits_{k\in K}(d_0^k,\varphi_k)=(0,0)$, it follows that
$\lim\limits_{k\in K}\tilde{d}^k=\lim\limits_{k\in
K}q^k=\lim\limits_{k\in K}d_1^k=0$.
\end{proof}

Now, we present the main result of this section.
\begin{theorem}\label{theorem1} Suppose that Assumptions \ref{assumption1},
\ref{assumption2} and \ref{assumption3} hold. Then Algorithm
\ref{algorithm1} either stops at a KKT point $x^k$ for problem (NCP)
after a finite number of iterations or generates an infinite
sequence $\{x^k\}$ of points such that each accumulation point $x^*$
(if it exists) of $\{x^k\}$ is a KKT point for problem (NCP),
furthermore, there exists an index set $K$ such that
$\{(x^k,\lambda^k)\}_K$ converges to the KKT pair $(x^*,\lambda^*)$
for problem (NCP).
\end{theorem}
\begin{proof} Choose an infinite index set $K$ such that
(\ref{21}) holds, and let matrix $R_k=(g_j(x^k),\ j\in J)$. In view
of $(x^k,d_0^k,\varphi_k)\rightarrow(x^*,0,0),\ k\in K$, we can
conclude that $J\subseteq I(x^*)=\{j\in I:\ f_j(x^*)=0\}$, and this
together with Assumption \ref{assumption2} shows that $R_k^TR_k$ is
nonsingular for $k\in K$ large enough, since
$R_k\stackrel{K}{\longrightarrow}R_*\triangleq(g_j(x^*),\ j\in J)$.
Again, from the KKT condition (\ref{4}), we have
$$g_0(x^k)+B_kd_0^k+R_k\lambda^k_{J}=0.$$ Thus, for $k\in K$
sufficiently large, we have
$$\lambda^k_{J}=-(R_k^TR_k)^{-1}R_k^T(g_0(x^k)
+B_kd_0^k)\rightarrow-(R_*^TR_*)^{-1}R_*^Tg_0(x^*)\stackrel{def}{=}\lambda^*_{J}.$$

If we denote the multiplier vector
$\lambda^*=(\lambda^*_{J},0_{I\setminus J})$, then
$\lim\limits_{k\in K}\lambda^k=\lambda^*.$ Furthermore, passing to
the limit $k\in K$ and $k\rightarrow\infty$ in (\ref{4}), it follows
that
$$g_0(x^*)+N_*\lambda^*=0,\ f_j(x^*)\leq 0,\ \lambda_j^*\geq 0,\ f_j(x^*)\lambda_j^*=0,\ j\in I,$$
which shows that $(x^*,\lambda^*)$ is a KKT pair for the problem
(NCP). \end{proof}
\begin{remark}\label{remark8} The global convergence Theorem \ref{theorem1} shows that
if the sequence $x^k$ generated by the proposed Algorithm
\ref{algorithm1} possesses a limit point $x^*$, then $x^*$ is a KKT
point for problem (NCP).
\end{remark}
\section{Superlinear convergence}\label{section4}
In this section, under some suitable assumptions, we first prove
that Algorithm \ref{algorithm1} possesses strong convergence, and
the iteration points always enter into the feasible set $X$ after a
finite number of iterations. Subsequently, the superlinear rate of
convergence is established without the strict complementarity. For
these purposes, we further make the following assumption.

\begin{assumption}\label{assumption4} (i)\ The functions $f_j(x)\
(j\in\{0\}\cup I)$ are all second-order
continuously differentiable.\\
(ii)\ The sequence $\{x^k\}$ generated by Algorithm \ref{algorithm1}
is bounded, and possesses an accumulation point $x^*$, such that the
KKT pair $(x^*,\lambda^*)$ satisfies the strong second-order
sufficient conditions (SSOSC for short), i.e.,
$$d^T\nabla^2_{xx}L(x^*,\lambda^*)d>0,\ \forall d\in R^n,\ d\neq 0,\ g_j(x^*)^Td=0,\ j\in I_*^+,$$
where $L(x,\lambda)=f_0(x)+\sum\limits_{j\in I}\lambda_jf_j(x),\
I_*^+=\{j\in I:\ \lambda_j^*>0\}.$
\end{assumption}
\begin{remark}\label{remark9} On the one hand, in some previously proposed SQP-type algorithms \cite{pt1987,zj2005},
to get the superlinear convergence of the proposed algorithm, one
has to ensure relation $J_k\equiv I(x^*)$ holds for $k$ large
enough. And the strict complementarity (i.e., $I_*^+=I(x^*)$) is an
very important condition for ensuring  $J_k\equiv I(x^*)$ holds. On
the other hand, it is well-known that the SSOSC is equivalent to the
second-order sufficient conditions (SOSC for short) under the
condition of the strict complementarity, however, this condition is
hard to verify in practise, and the positive space (i.e., the
critical directions set) of the SOSC is smaller than the SSOSC's. In
this paper, in order to obtain the strong and superlinear
convergence of Algorithm \ref{algorithm1}, combining a slight strong
assumption SSOSC, we can still avoid the Maratos effect under
$I_*^+\subseteq J_k\subseteq I(x^*)$ (see Lemma \ref{lemma8},
Theorem \ref{theorem3}).
\end{remark}

First, under the stated assumptions, we have the following theorem.
\begin{theorem}\label{theorem2}Suppose that Assumptions \ref{assumption2}, \ref{assumption3}
and \ref{assumption4} hold.
Then\\
(i)
$\lim\limits_{k\rightarrow+\infty}d_0^k=\lim\limits_{k\rightarrow+\infty}d_1^k=
\lim\limits_{k\rightarrow+\infty}\tilde{d}^k=
\lim\limits_{k\rightarrow+\infty}q^k=\lim\limits_{k\rightarrow+\infty}d^k=0$
and
$\lim\limits_{k\rightarrow+\infty}||x^{k+1}-x^k||=0$;\\
(ii)$\lim\limits_{k\rightarrow+\infty} x^k=x^*,$ i.e., Algorithm
\ref{algorithm1}
is said to be strongly convergent in this sense; \\
(iii) $\lim\limits_{k\rightarrow+\infty}\varphi_k=\varphi(x^*)=0$,
$\lim\limits_{k\rightarrow+\infty}\lambda^k=\lambda^*$.
\end{theorem}
\begin{proof} (i)\ Since $\{x^k\}$ is bounded, from Lemma
\ref{lemma7}, it follows that any subsequence of $\{(d_0^k, d_1^k,
\tilde{d}^k, q^k)\}$ must possesses an accumulation point
$(0,0,0,0)\in R^{4n}$, which implies
$\lim\limits_{k\rightarrow+\infty}(d_0^k,d_1^k,\tilde{d}^k,q^k)$
$=(0,0,0,0)$. Furthermore, from the mechanism of Algorithm
\ref{algorithm1}, we obtain
$$\lim\limits_{k\rightarrow+\infty}d^k=0,\ \
\lim\limits_{k\rightarrow+\infty}||x^{k+1}-x^k||=\lim\limits_{k\rightarrow+\infty}||t_kd^k||=0.$$

(ii) First, under Assumptions \ref{assumption2} and
\ref{assumption4}, by Proposition 4.1 in \cite{pt1987}, it is known
that the given limit point $x^*$ is an isolated KKT point of problem
(NCP). Therefore, we know from Theorem \ref{theorem1} that $x^*$ is
an isolated limit point of $\{x^k\}$, this combining
$||x^{k+1}-x^k||\rightarrow 0$, it follows that $x^k\rightarrow x^*$
(The details can be found in \cite{j2000}).

(iii) In view of Theorem \ref{theorem1}, we know that the given
accumulation point $x^*$ is a KKT point for problem (NCP), thus it
follows $\varphi(x^*)=0$. Furthermore, combining the monotone
property and boundedness of $\{\varphi_k\}$, we have
$\lim\limits_{k\rightarrow+\infty}\varphi_k=\varphi(x^*)=0.$
Moreover, from the proof of Theorem \ref{theorem1} and parts (i) and
(ii), one can conclude that each accumulation point of sequence
$\{\lambda_k\}$ is a KKT multiplier associated with $x^*$, this
together with the uniqueness of the KKT multiplier implies that
$\lim\limits_{k\rightarrow+\infty}\lambda^k=\lambda^*$.
\end{proof}

The relations established in the following lemma are very important
in the subsequent discussion.
\begin{lemma}\label{lemma8}
Suppose that Assumptions \ref{assumption2}, \ref{assumption3} and \ref{assumption4} all hold. Then\\
(i) $||d_1^k||=O(||d_0^k||^2)+O(\varphi_k^\sigma)$,
$||d_1^k||^2=O(||d_0^k||^4)
+o(\varphi_k^\sigma)$, $||h_1^k||=O(||d_0^k||^2)+O(\varphi_k^\sigma)$;\\
(ii) $I_*^+\subseteq J_k\subseteq I(x^*)$ for $k$
large enough;\\
(iii) the relationship (\ref{13}) is satisfied for $k$ large enough.
\end{lemma}
\begin{proof} (i) In view of $\tilde{F}(x^k+d_0^k)=O(||d_0^k||^2)$,
the proof is elementary from (\ref{11}) and Lemma \ref{lemma6}(ii)
as well as Theorem \ref{theorem2}.

(ii) For $j\not\in I(x^*)$, i.e., $f_j(x^*)<0$. Since
$\lim\limits_{k\rightarrow+\infty}(x^k, d_0^k)=(x^*, 0)$, there
exists a constant $\bar{\xi}>0$ such that
$\bar{f}_j(x^k)=f_j(x^k)\leq-\bar{\xi}<0$ holds for $k$ large
enough. It then follows that $\bar{f}_j(x^k)+g_j(x^k)^Td_0^k\leq
-\frac{1}{2}\bar{\xi}<0$ for $k$ large enough, which shows that
$j\not\in J_k$, thus $J_k\subseteq I(x^*)$. Furthermore, it follows
from Theorem \ref{theorem2} that
$\lim\limits_{k\rightarrow+\infty}\lambda_{I_*^+}^k=\lambda_{I_*^+}^*>0$,
so, $\lambda_{I_*^+}^k>0$ and $I_*^+\subseteq J_k$ hold for $k$
sufficiently large.

(iii) From (\ref{5}) and Assumption \ref{assumption3}, we only need
to show that
$$-\frac{a}{2}||d_0^k||^2\leq\zeta\min\{-||d_0^k||^\delta,\
-||d^k||^\delta\}+\xi\varphi_k^\varrho$$ holds for $k$ large enough.

First, in view of $\delta>2$ and Lemma \ref{lemma7}, it follows that
\begin{equation}\label{24}-\frac{a}{2}||d_0^k||^2\leq-\zeta||d_0^k||^\delta\leq-\zeta||d_0^k||^\delta
+\xi\varphi_k^\varrho,\end{equation} for $k$ sufficiently large and
$\zeta,a>0$.

Second, in view of $||d^k||=||d_0^k+d_1^k||\leq||d_0^k||+||d_1^k||$,
$\delta>2$, part (i) and Lemma \ref{lemma7}, we obtain
$$||d^k||^\delta\leq(||d_0^k||+O(||d_0^k||^2)+O(\varphi_k^\sigma))^\delta=
||d_0^k||^\delta+o(||d_0^k||^2)+o(\varphi_k^\sigma)$$ for $k$ large
enough. This together with $\delta>2$ and $\varrho<\sigma$ implies
that
$$\begin{array}{ll}-\zeta||d^k||^\delta+\xi\varphi_k^\varrho+\frac{a}{2}||d_0^k||^2&\geq-
\zeta||d_0^k||^\delta+\frac{a}{2}||d_0^k||^2+o(||d_0^k||^2)+\xi\varphi_k^\varrho+o(\varphi_k^\sigma)\\
&=\frac{a}{2}||d_0^k||^2+o(||d_0^k||^2)+\xi\varphi_k^\varrho+o(\varphi_k^\sigma)\\
&\geq 0\end{array}$$ holds for $k$ large enough, i.e.,
\begin{equation}\label{25}-\frac{a}{2}||d_0^k||^2\leq-\zeta||d^k||^\delta+\xi\varphi_k^\varrho,\
\ {\rm{for}}\ k\ \rm{large\ enough}.\end{equation} Finally,
combining (\ref{24}), (\ref{25}) and (\ref{5}), the conclusion (iii)
holds for $k$ large enough. \end{proof}

To ensure that the steplength $t_k\equiv1$ for $k$ large enough
without the strict complementary assumption, an additional
assumption as follows is necessary.

\begin{assumption}\label{assumption5} Suppose that the KKT pair
$(x^*,\lambda^*)$ and the matrix $B_k$ satisfy
$$||(\nabla^2_{xx}L(x^*,\lambda^*)-B_k)d_0^k||=o(||d_0^k||).$$
\end{assumption}

\begin{remark}\label{remark10}According to Theorem \ref{theorem5}, it holds that
$${\bf{Assumption\ 5}}\Longleftrightarrow||(\nabla^2_{xx}L(x^k,\lambda^k)-B_k)d_0^k||=o(||d_0^k||).$$
\end{remark}

\begin{theorem}\label{theorem3} Suppose that Assumptions \ref{assumption2}, \ref{assumption3},
\ref{assumption4} and \ref{assumption5} are all satisfied. Then, the
inequality (\ref{14}) in Step 3 is always satisfied for $t = 1$ and
sufficiently large $k$. Therefore, Step 4 and Step 5 are no longer
performed in Algorithm \ref{algorithm1}, and Algorithm
\ref{algorithm1} always performs cycle I.
\end{theorem}
\begin{proof} First of all, in view of Lemma \ref{lemma8}(iii), it is
sufficient to prove that the inequality (\ref{14}) holds for $t = 1$
and sufficiently large $k$.

Discuss the second group of inequalities and the last group of
inequalities of (\ref{14}).

For $j\not\in I(x^*)$, i.e., $f_j(x^*)<0$, in view of
$(x^k,d_0^k,d_1^k,\varphi_k)\rightarrow(x^*,0,0,0)\
(k\rightarrow+\infty)$, we obtain $d^k=d_0^k+d_1^k\rightarrow0\
(k\rightarrow+\infty)$, thus, we conclude that the second group of
inequalities and the last group of inequalities of (\ref{14}) are
both satisfied for $t=1$ and $k$ large enough.

For $j\in I(x^*)$, since
$\lim\limits_{k\rightarrow+\infty}f_j(x^k)=f_j(x^*)=0$ and
$\lim\limits_{k\rightarrow+\infty}\varphi_k=0$ as well as (\ref{3})
and (\ref{7}), it follows that
$$\bar{f}_j(x^k)\rightarrow 0,\ D_j^k\rightarrow 0,\
D_j^k=o(|\bar{f}_j(x^k)+g_j(x^k)^Td_0^k|)+o(||d_0^k||).$$
 On the
other hand, from (\ref{11}) and Lemma \ref{lemma8}(i), we have
\begin{equation}\label{26}\begin{array}{ll}g_j(x^k)^Td_1^k&=-||d_0^k||^\tau-\varphi_k^\sigma-f_j(x^k+d_0^k)
+f_j(x^k)+g_j(x^k)^Td_0^k+D_j^kh_{1j}^k\\
&=-||d_0^k||^\tau-\varphi_k^\sigma-f_j(x^k+d_0^k)+f_j(x^k)+g_j(x^k)^Td_0^k\\
&\ \ \
+o(|\bar{f}_j(x^k)+g_j(x^k)^Td_0^k|)+O(||d_0^k||^3)+o(\varphi_k^\sigma).\end{array}\end{equation}
Then, from Taylor expansion, Lemma \ref{lemma8}(i), (\ref{26}) and
$\tau\in(2,3)$, we obtain
\begin{equation}\label{27}\begin{array}{ll}f_j(x^k+d^k)&=f_j(x^k+d_0^k)+g_j(x^k+d_0^k)^Td_1^k+O(||d_1^k||^2)\\
&=f_j(x^k+d_0^k)+g_j(x^k)^Td_1^k+O(||d_0^k||^3)+o(\varphi_k^\sigma)\\
&=-||d_0^k||^\tau-\varphi_k^\sigma+f_j(x^k)+g_j(x^k)^Td_0^k+o(|\bar{f}_j(x^k)+g_j(x^k)^Td_0^k|)\\
&\ \ \ +O(||d_0^k||^3)+o(\varphi_k^\sigma),\end{array}\end{equation}
which further implies that
\begin{equation}\label{28}f_j(x^k+d^k)=\left\{\begin{array}{ll}-||d_0^k||^\tau-\varphi_k^\sigma-|\bar{f}_j(x^k)
+g_j(x^k)^Td_0^k|+o(|\bar{f}_j(x^k)+g_j(x^k)^Td_0^k|)\\
\ \ \ \ +O(||d_0^k||^3)+o(\varphi_k^\sigma),\ \ \ \ \ \ {\rm{if}}\
j\in
I(x^*)\cap I^-;\\
-||d_0^k||^\tau-\varphi_k^\sigma+\varphi_k-|\bar{f}_j(x^k)+g_j(x^k)^Td_0^k|+o(|\bar{f}_j(x^k)+g_j(x^k)^Td_0^k|)\\
\ \ \ \ +O(||d_0^k||^3)+o(\varphi_k^\sigma),\ \ \ \ \ \ {\rm{if}}\
j\in I(x^*)\cap I^+.\end{array}\right.\end{equation} This shows that
$f_j(x^k+d^k)\leq 0$, for $j\in I(x^*)\cap I^-$ and $k$ large
enough.

For $j\in I(x^*)\cap I^+$, from (\ref{28}) and $\tau\in(2,3)$ as
well as $\alpha\in(0,\frac{1}{2})$, it follows that
\begin{equation}\label{29}\begin{array}{ll}f_j(x^k&+d^k)-\varphi_k+\alpha(||d_0^k||^\tau+\varphi_k^\sigma)\\
&=-(1-\alpha)(||d_0^k||^\tau+\varphi_k^\sigma)-|\bar{f}_j(x^k)+g_j(x^k)^Td_0^k|\\
&\ \ \ \ +o(|\bar{f}_j(x^k)+g_j(x^k)^Td_0^k|)+o(\varphi_k^\sigma)+O(||d_0^k||^3)\\
&\leq 0.\end{array}\end{equation} Thus, summarizing the analysis
above, we have proved that the inequality of (\ref{14}) except the
first one is satisfied for $t=1$ and $k$ large enough.

Finally, we will show that the first inequality of (\ref{14}) holds
for $t=1$ and $k$ large enough. From Taylor expansion and Lemma
\ref{lemma8}(i), we have
\begin{equation}\label{30}\begin{array}{ll}\Delta_k&\triangleq
f_0(x^k+d^k)-f_0(x^k)-\alpha
g_0(x^k)^Td_0^k-\rho(1-\alpha)\varphi_k^\theta\\
&=g_0(x^k)^Td^k+\frac{1}{2}(d^k)^T\nabla_{xx}^2f_0(x^k)d^k-\alpha
g_0(x^k)^Td_0^k-\rho(1-\alpha)\varphi_k^\theta+o(||d^k||^2)\\
&=g_0(x^k)^T(d_0^k+d_1^k)+\frac{1}{2}(d_0^k)^T\nabla_{xx}^2f_0(x^k)d_0^k-\alpha
g_0(x^k)^Td_0^k-\rho(1-\alpha)\varphi_k^\theta\\&\ \
+o(||d_0^k||^2)+o(\varphi_k^\sigma).\end{array}\end{equation}

Then, from the KKT conditions (\ref{4}) and Lemma \ref{lemma8}(i),
we have

\begin{equation}\label{31}g_0(x^k)^Td_0^k=-(d_0^k)^TB_kd_0^k-\sum\limits_{j\in
J_k}\lambda_j^kg_j(x^k)^Td_0^k,\end{equation}

\begin{equation}\label{32}\begin{array}{ll}
g_0(x^k)^T(d_0^k+d_1^k)=-(d_0^k)^TB_kd_0^k-\sum\limits_{j\in
J_k}\lambda_j^kg_j(x^k)^T(d_0^k+d_1^k)+o(||d_0^k||^2)+o(\varphi_k^\sigma).\end{array}\end{equation}

For $j\in J_k\subseteq I(x^*)$, it follows that
$\bar{f}_j(x^k)+\nabla f_j(x^k)^Td_0^k=0$, which together with
(\ref{28}) and $\varphi_k=o(\varphi_k^\sigma)$ implies
\begin{equation}\label{33}\begin{array}{ll}
f_j(x^k+d^k)=-||d_0^k||^\tau-\varphi_k^\sigma+o(||d_0^k||^2)+o(\varphi_k^\sigma),\
j\in J_k.\end{array}\end{equation}

Again, from Taylor expansion and Lemma \ref{lemma8}(i) as well as
$\varphi_k=o(\varphi_k^\sigma)$, we get
\begin{equation}\label{34}\begin{array}{ll}f_j(x^k+d^k)&=f_j(x^k)+g_j(x^k)^T(d_0^k+d_1^k)+
\frac{1}{2}(d_0^k)^T\nabla_{xx}^2f_j(x^k)d_0^k+o(||d_0^k||^2)+o(\varphi_k^\sigma)\\
&=\bar{f}_j(x^k)+g_j(x^k)^T(d_0^k+d_1^k)+\frac{1}{2}(d_0^k)^T\nabla_{xx}^2f_j(x^k)d_0^k+o(||d_0^k||^2)
+o(\varphi_k^\sigma).\end{array}\end{equation}

Combining (\ref{33}) and (\ref{34}), it follows that
\begin{equation}\label{35}\begin{array}{ll}-\sum\limits_{j\in
J_k}\lambda_j^kg_j(x^k)^T(d_0^k+d_1^k)&=\sum\limits_{j\in
J_k}\lambda_j^k\bar{f}_j(x^k)+\frac{1}{2}\sum\limits_{j\in
J_k}\lambda_j^k(d_0^k)^T\nabla_{xx}^2f_j(x^k)d_0^k\\
&\ \ \ \
+o(||d_0^k||^2)+O(\varphi_k^\sigma).\end{array}\end{equation}

Thus, substituting (\ref{35}) into (\ref{32}), we have
\begin{equation}\label{36}\begin{array}{ll}g_0(x^k)^T(d_0^k+d_1^k)=&-(d_0^k)^TB_kd_0^k+\sum\limits_{j\in
J_k}\lambda_j^k\bar{f}_j(x^k)+\frac{1}{2}(d_0^k)^T(\nabla^2_{xx}L(x^k,\lambda_k)\\&-\nabla^2_{xx}f_0(x^k))d_0^k+o(||d_0^k||^2)+O(\varphi_k^\sigma).\end{array}\end{equation}

In addition, in view of $g_j(x^k)^Td_0^k=-\bar{f}_j(x^k)$, $j\in
J_k$, from (\ref{31}), it follows that

\begin{equation}\label{37}g_0(x^k)^Td_0^k=-(d_0^k)^TB_kd_0^k+\sum\limits_{j\in
J_k}\lambda_j^k\bar{f}_j(x^k).\end{equation}

Now, substituting (\ref{36}) and (\ref{37}) into (\ref{30}), we have
$$\begin{array}{ll}\Delta_k&=-(d_0^k)^TB_kd_0^k+\sum\limits_{j\in J_k}\lambda_j^k\bar{f}_j(x^k)
+\frac{1}{2}(d_0^k)^T\nabla_{xx}^2L(x^k,\lambda_k)d_0^k-\alpha g_0(x^k)^Td_0^k\\
&\ \ \ \ -\rho(1-\alpha)\varphi_k^\theta+o(||d_0^k||^2)+O(\varphi_k^\sigma)\\
&=(\alpha-\frac{1}{2})(d_0^k)^TB_kd_0^k+(1-\alpha)\sum\limits_{j\in
J_k}\lambda_j^k\bar{f}_j(x^k)+\frac{1}{2}(d_0^k)^T(\nabla_{xx}^2L(x^k,\lambda_k)\\
&\ \ \ \
-B_k)d_0^k-\rho(1-\alpha)\varphi_k^\theta+o(||d_0^k||^2)+O(\varphi_k^\sigma).\end{array}$$

This together with Assumptions \ref{assumption3}, \ref{assumption5},
$\lambda_j^k\bar{f}_j(x^k)\leq 0$ and $\alpha\in(0,\frac{1}{2})$ as
well as $\theta<\sigma$ shows that
$$\Delta_k\leq(\alpha-\frac{1}{2})a||d_0^k||^2+o(||d_0^k||^2)-\rho(1-\alpha)\varphi_k^\theta
+o(\varphi_k^\theta)\leq 0$$ holds for $k$ large enough. Hence, the
first inequality of (\ref{14}) holds for $t=1$ and $k$ large enough.
The whole proof is completed.\end{proof}

\begin{remark}\label{remark11} In order to overcome the Maratos
effect, under some suitable assumptions, some norm relations of
directions and a weaker set relation corresponding to the strict
complementarity (i.e., $I_*^+=I(x^*)$) are given in Lemma
\ref{lemma8} firstly. Then, by these results, we obtain the Theorem
\ref{theorem3}, i.e., the inequality (\ref{14}) is always satisfied
for $t=1$ and $k$ large enough (i.e., very close to the solution of
the problem). So, the Maratos effect can be overcame in our paper.
\end{remark}

According to Theorem \ref{theorem3} and its proof for case of
$j\not\in I(x^*)$ as well as relationship (\ref{28}), the following
lemma holds immediately.
\begin{theorem}\label{theorem4} Under the assumptions stated in Theorem \ref{theorem3},
we have $\varphi_{k+1}\equiv0$ after a finite number of iterations,
i.e., $x^{k+1}\in X$ for $k$ large enough.
\end{theorem}

At the end of this section, based on Theorems \ref{theorem3} and
\ref{theorem4} as well as Lemma \ref{lemma8}, using Theorem 3.1.3 in
\cite{j2000}, we have the superlinear convergence of Algorithm
\ref{algorithm1} immediately as follows.

\begin{theorem}\label{theorem5} Suppose that Assumptions
\ref{assumption2}, \ref{assumption3}, \ref{assumption4} and
\ref{assumption5} are all satisfied. Then,
$||x^{k+1}-x^*||=o(||x^k-x^*||)$, i.e., Algorithm \ref{algorithm1}
is superlinearly convergent. \end{theorem}
\section{Numerical experiments}\label{section5}
In this section, in order to illustrate the computational efficiency
of Algorithm \ref{algorithm1}, some preliminary numerical results
are reported, and the computing results show that Algorithm
\ref{algorithm1} is effective. The algorithm was implemented by
using Matlab 7.5 on Windows XP platform, and on a PC with 1.99 GHZ
CPU. The approximation Hession matrix $B_k$ is updated by the BFGS
formula described in \cite{pm1991}.
%

During the numerical experiments, the parameters are selected as
follows: $$\left\{\begin{array}{ll}\gamma=\eta=0.5,\
\theta=\varrho=0.4,\ \sigma=0.6,\ \ \xi=1,\ \zeta=0.2,\\
\alpha=0.3,\ \rho=1.5,\ \delta=3,\ \tau=2.5,\
\varepsilon=0.5^{3}.\end{array}\right.$$ We test some problems which
are taken from \cite{hs1980,s1987}. In addition, we further test
Svanberg problems in different dimensions and with different initial
points, which are taken from \cite{got2003}. Execution is terminated
if the norm of $d_0^k$ is less than a given constant $\epsilon>0$
and $\varphi_k=0$. The columns of the following tables have the
following meanings:

Prob: the number of the test problem in \cite{hs1980,s1987};

$n/m$: the number of variables/inequality constraints of the
problem;

Code: the name of the algorithm;

NF0: the number of objective function evaluations;

NIO/NII: the number of iterations out of/within the feasible set;

NI: the total number of iterations, i.e., NI=NIO+NII;

NF: the number of all constraint functions evaluations;

FV: the objective function value at the final iteration point;

CPU: the CPU time (second).

Finally, an ``$-$" in the following tables indicates that the
corresponding information is not given in the corresponding
references.

\textbf{Experiment 1} ({\it for small-scale problems}). For this
part, in order to show the computational efficiency of Algorithm
\ref{algorithm1} (denoted by ALG \ref{algorithm1}), we test some
small-scale problems and compare ALG \ref{algorithm1} with some
other algorithms, and the numerical results are given in Tables 1-4.

In Tables 1 and 2, ALG \ref{algorithm1} is compared with ALGO
\cite{jzth2006} and SNQP \cite{jkzt2009} for the same test problems,
the stopping criterion threshold $\epsilon$ and initial iteration
points are the same as that reported in \cite{jzth2006} and
\cite{jkzt2009}, respectively. From the viewpoint of the numbers of
NIO, it follows that ALG \ref{algorithm1} can always enter into $X$
after relatively small iterations. Furthermore, from the viewpoint
of the numbers of NIO, NII and FV, the results show that ALG
\ref{algorithm1} is obviously better than ALGO for most of test
problems. The performance of ALG \ref{algorithm1} in terms of NII is
better than SNQP except problems 33 and 76.

Tables 3 and 4 gives the compared numerical results for ALG
\ref{algorithm1} and ALG 3.1 as well as ALG 3.2 \cite{ph1991}. The
test problems and stopping criterion threshold are the same as in
\cite{ph1991}. The numerical results in Tables 3 and 4 show that ALG
\ref{algorithm1} can always enter into $X$ after small iterations,
and ALG \ref{algorithm1} is more better than ALG 3.1 and ALG 3.2 for
the test problems.
\begin{center}
{\emph{{\rm Table 1.} Numerical results for Experiment 1-I}}\\
{\scriptsize
\begin{tabular}{ c c c c c c c c c c c c c c c c }
 \hline\hline
   Prob   &$n/m$ &Initial point     &Code    &NIO      &NII    &NF0       &NF      &FV       &CPU\\
 \hline
   012     &2/1   &$(6,6)^T$ &ALG \ref{algorithm1}  &17   &3  &21     &41      &$-3.0000000E+01$   &$0.06$\\
             &      &                &ALGO     &25       &28     &29        &57      &$-3.0000000E+01$   &$-$\\
             &      &                &SNQP     &7       &12      &12         &29     &$-2.9999999E+01$ &$-$\\
 \hline
   029       &3/1   &$(-4,-4,-4)^T$&ALG \ref{algorithm1} &3 &9  &13    &46     &$-2.2627417E+01$   &$0.05$\\
             &      &              &ALGO  &1 &11 &14 &27 &$-2.2627417E+01$ &$-$\\
             &      &                 &SNQP    &1      &12    &17         &42      &$-2.2627416E+01$     &$-$\\
 \hline
   031        &3/7   &$(2,4,7)^T$&ALG \ref{algorithm1} &1 &16     &18     &309  &$6.0000000E+00$   &$0.06$\\
             &     &　　　　　　　　　&ALGO    &4      &20     &23     &43    &$6.0000000E+00$   &$-$\\
             &     &                  &SNQP    &1      &19     &20     &52    &$6.0000089E+00$  &$-$\\
 \hline
   033        &3/6   &$(2,4,6)^T$&ALG \ref{algorithm1} &1   &9     &11    &134  &$-4.5857864E+00$   &$0.05$\\
   　　　　　&       &               &ALGO &2       &16      &17   &67   &$-4.5857863E+00$           &$-$\\
             &       &$(1,4,6)^T$&ALG \ref{algorithm1} &1  &44     &46    &570    &$-4.5857864E+00$ &$0.33$\\
             &       &               &SNQP  &2   &21     &23     &116   &$-4.5857290E+00$  &$-$\\
 \hline
  034        &3/8   &$(2,2,2)^T$&ALG \ref{algorithm1} &5 &10 &16     &166       &$-8.3403245E-01$   &$0.06$\\
              &        &　　　　　&ALGO    &8      &26     &68    &198       &$-8.3403244E-01$           &$-$\\
 \hline
   035        &3/4   &$(1,2,3)^T$&ALG \ref{algorithm1} &1  &6   &8 &67          &$1.1111111E-01$   &$0.03$\\
              &      &           &ALGO  &8 &11 &12 &0            &$-3.4500000E+00$ &$-$\\
              &     &              &SNQP  &4     &9   &9   &0    &$1.1111111E-01$ &$-$\\
 \hline
   043        &4/3   &$(-10,2,-8,5)^T$&ALG \ref{algorithm1} &9  &5  &15  &95    &$-4.4000000E+01$   &$0.06$\\
             &          &　　　　　　　&ALGO    &23      &26    &27     &163      &$-4.4000000E+01$     &$-$\\
             &      &$(0,2,2,4)^T$&ALG \ref{algorithm1} &7   &9    &17    &135   &$-4.4000000E+01$  &$0.08$\\
             &      &                  &SNQP   &1   &11   &11   &69   &$-4.3999999E+01$    &$-$\\
 \hline
   044        &4/10  &$(-20,-20,
                      -20,-20)^T$&ALG \ref{algorithm1}  &4  &10  &15   &296  &$-1.5000000E+01$   &$0.08$\\
   　　　　　&      &                   &ALGO    &7    &15    &18     &0   &$-1.5000000E+01$           &$-$\\
 \hline
   066        &3/8   &$(0,0,100)^T$&ALG \ref{algorithm1} &10  &54  &65  &1067  &$5.1816327E-01$   &$0.48$\\
             &         &              &ALGO  &34      &39     &40     &161   &$5.1816327E-01$           &$-$\\
 \hline
   076        &4/7   &$(1,2,3,4)^T$&ALG \ref{algorithm1} &5  &16  &22     &345  &$-4.6818182E+00$   &$0.11$\\
            &          &　　　　　　　&ALGO   &6       &14     &15    &0  &$-4.6818182E+00$   &$-$\\
            &       &                 &SNQP   &2       &14      &14    &0  &$-4.6818171E+00$ &$-$\\
 \hline
   100       &7/4   &$\begin{array}{ll}(0,3,-3,\\
                     3,0,1,0)^T\end{array}$&ALG \ref{algorithm1} &18 &39  &58      &861  &$6.8256637E+02$   &$0.61$\\
   　　　　　&      &                  &ALGO  &4       &33      &47    &363  &$6.8063006E+02$           &$-$\\
 \hline\hline
\end{tabular}
}
\end{center}
\newpage
\begin{center}
{\emph{{\rm Table 2.} Numerical results for Experiment 1-I-continued}}\\
{\scriptsize
\begin{tabular}{ c c c c c c c c c c c c c c c c }
 \hline\hline
   Prob   &$n/m$ &Initial point     &Code    &NIO      &NII    &NF0       &NF      &FV       &CPU\\
\hline
  113        &10/8  &$\begin{array}{ll}(4,10,10,2,0,\\
                     11,4,0,12,10)^T\end{array}$&ALG \ref{algorithm1} &12 &4  &17 &378  &$2.4306209E+01$   &$0.14$\\
             &      &                       &ALGO  &6     &17    &21   &205  &$2.4306209E+01$           &$-$\\
             &      &$\begin{array}{ll}(0,2,9,5,0,\\
                     1,9,8,-10,10)^T\end{array}$&ALG \ref{algorithm1} &9  &7  &17  &355   &$2.4306585E+01$  &$0.17$\\
             &       &                     &SNQP    &9    &29     &29   &258     &$2.4306211E+01$ &$-$\\
 \hline
   264       &4/3    &$(8,-5,6,-4)^T$&ALG \ref{algorithm1} &18  &5  &24  &142  &$-4.3987578E+01$   &$0.13$\\
           &         &                 &ALGO   &19     &26     &27    &161   &$-4.3999999E+01$    &$-$\\
           &         &$(0,0,0,10)^T$&ALG \ref{algorithm1}  &17   &5  &23  &146  &$-4.3987578E+01$ &$0.16$\\
           &         &                 &SNQP    &4      &16   &17    &122   &$-4.4113405E+01$  &$-$\\
 \hline\hline
\end{tabular}
}
\end{center}

\begin{center}
\emph{{\rm Table 3.} Numerical results for Experiment 1-II}\\
{\scriptsize
\begin{tabular}{ c c c c c c c c c c c c c c c c }
\hline\hline
Prob           &$n/m$   &Initial point       &Code          &NIO+NII  &NF0+NF &NDF0+NG  &CPU\\
\hline
Rosen/Suzuki-1 &4/3     &$(0,0,0,0)^T$&ALG A      &0+17    &18+117       &18+54     &0.09\\
               &        &                    &ALG 3.1  &76  &1185   &308      &2.88\\
               &        &                    &ALG 3.2  &77  &1225   &312      &3.04\\
Rosen/Suzuki-2 &4/3     &$(2,4,8,1)^T$&ALG A  &9+10    &20+120       &20+60  &0.09\\
               &        &                    &ALG 3.1  &68  &1034   &276      &2.36\\
               &        &                    &ALG 3.2  &55  &793    &224      &2.06\\
\hline
Wong Problem-1 &7/4     &$(1,2,0,4,0,1,1)^T$&ALG A  &0+24    &25+228       &25+100 &0.16\\
               &        &                    &ALG 3.1  &157  &17756   &790      &21.48\\
               &        &                    &ALG 3.2  &157  &17759   &790      &21.54\\
Wong Problem-2 &7/4     &$(3,3,0,5,1,3,0)^T$&ALG A &9+48  &58+624       &58+232        &0.38\\
               &        &                    &ALG 3.1  &171  &18677   &860      &22.48\\
               &        &                    &ALG 3.2  &151  &16921    &760     &20.88\\
\hline
Quadratic Problem-1 &2/2     &$(-0.3,0.0)^T$&ALG A  &0+7    &8+38       &8+16        &0.03\\
               &        &                    &ALG 3.1  &48  &292   &147      &0.76\\
               &        &                    &ALG 3.2  &49  &301   &150      &0.76\\
Quadratic Problem-2 &2/2     &$(2.2,1.6)^T$&ALG A &6+4    &11+42       &11+22        &0.06\\
               &        &                    &ALG 3.1  &50  &314   &153      &0.72\\
               &        &                    &ALG 3.2  &43  &286    &132     &0.68\\
\hline\hline
\end{tabular}
}
\end{center}
\newpage
\begin{center}
\emph{{\rm Table 4.} Numerical results for Experiment 1-II-continued}\\
{\scriptsize
\begin{tabular}{ c c c c c c c c c c c c c c c c }
\hline\hline
Prob           &$n/m$   &Initial point       &Code          &NIO+NII  &NF0+NF &NDF0+NG  &CPU\\
\hline
TFI1 Problem-1 &3/1     &$(-10,0,0)^T$ &ALG A  &0+12    &13+32       &13+13        &0.05\\
               &        &                    &ALG 3.1  &35  &2663   &792      &1.46\\
               &        &                    &ALG 3.2  &35  &2663   &792      &1.46\\
TFI1 Problem-2 &3/1     &$(1,1,1)^T$&ALG A &1+7 &9+17      &9+9        &0.05\\
               &        &                    &ALG 3.1  &23  &2301   &528      &0.98\\
               &        &                    &ALG 3.2  &22  &2289   &506     &0.96\\
\hline
TFI2 Problem-1 &3/1     &$(2,2,2)^T$&ALG A    &0+6    &7+23       &7+9        &0.04         \\
               &        &                    &ALG 3.1  &30  &1343   &682      &1.04\\
               &        &                    &ALG 3.2  &30  &1343   &682      &1.04\\
TFI2 Problem-2 &3/1     &$(0,0,0)^T$&ALG A  &4+7    &12+25       &12+18        &0.06        \\
               &        &                    &ALG 3.1  &48  &2129  &1078      &1.58\\
               &        &                    &ALG 3.2  &48  &2135   &1078     &1.60\\
\hline\hline
\end{tabular}
}
\end{center}

\textbf{Experiment 2} ({\it for middle-large-scale problems}).
Considering the all of tested problems above are all relatively
small, we further test the Svanberg problems \cite{got2003}
problems, some of them are larger and therefore interesting. The
experiment results are given in Tables 5-6.

In Table 5, the performance of ALG \ref{algorithm1} is compared with
SNQP, ALGO, FSLE \cite{ylq2003}. The initial iteration
 points and the stopping criterion threshold are the same as that
  reported in \cite{ylq2003}. From the results in Table 5, in
viewpoint of NII and NF0, it follows that algorithm ALG
\ref{algorithm1} performs better than FSLE, SNQP and ALGO in most
cases for problems Svanberg.

In Table 6, we further test Svanberg problems (in different
dimensions) for some infeasible initial points, and the stopping
criterion threshold is $\epsilon= 10^{-6}$ and $\varphi_k=0$. The
results show that our algorithm ALG \ref{algorithm1} is always
successful for all cases, and the iteration points can enter into
the feasible set so faster. In view of NIO, NII and CPU, it follows
that our algorithm is effective.
\newpage
\begin{center}
{\emph{{\rm Table 5.} Numerical results for Experiment 2-I}\\
\scriptsize
\begin{tabular}{ c  c  c  c  c  c  c  c  c  c  c  c  c  c  c  c }
 \hline\hline
   Prob      &$n/m$   &Initial point         &Code     &NI0       &NII        &NF0       &NF        &FV  &CPU\\
\hline
 Svanberg-10   &10/30   &$(0,0,\ldots,0)^T$&ALG \ref{algorithm1}&0  &16  &17  &1140     &15.731517 &$0.34$\\
               &        &                    &SNQP   &0      &28         &28        &1753    &15.731533  &$-$\\
               &        &                   &ALGO    &0      &15         &21        &1050     &15.731517 &$-$\\
               &        &                  &FSLE     &0      &36         &227       &258      &15.731517 &$-$\\
\hline
Svanberg-30    &30/90   &$(0,0,\ldots,0)^T$&ALG \ref{algorithm1}&0 &25  &26  &5490       &49.142526&$1.77$\\
               &        &                    &SNQP   &0     &27         &27        &4975       &49.142545&$-$\\
               &        &                      &ALGO  &0    &26         &38        &5670       &49.142526&$-$\\
               &        &                     &FSLE   &0    &101        &777       &864        &49.142526&$-$\\
\hline
 Svanberg-50   &50/150  &$(0,0,\ldots,0)^T$&ALG \ref{algorithm1}&0  &33  &34  &11550     &82.581912&$5.91$\\
              &&                              &SNQP &0     &37         &37        &11762       &82.581928&$-$\\
               &        &                     &ALGO &0     &35         &51        &12750       &82.581912&$-$\\
               &        &                    &FSLE   &0    &108        &881       &968         &82.581912&$-$\\
 \hline
 Svanberg-80   &80/240  &$(0,0,\ldots,0)^T$&ALG \ref{algorithm1}&0  &42 &43  &24720    &132.749819&$15.38$\\
              &&                             &SNQP &0     &47         &47        &24100       &132.749830&$-$\\
               &        &                    &ALGO  &0    &47         &68        &27360       &132.749819&$-$\\
               &        &                    &FSLE  &0    &190        &1666      &1835        &132.749819&$-$\\
 \hline
 Svanberg-100  &100/300  &$(0,0,\ldots,0)^T$&ALG \ref{algorithm1}&0  &46 &91 &53700    &166.197172&$26.38$\\
               &         &                   &SNQP &0  &46         &46        &27880       &166.197199&$-$\\
               &         &                &ALGO  &0      &53         &66        &35400       &166.197171&$-$\\
               &         &               &FSLE   &0      &178        &1628      &1782        &166.197171&$-$\\
 \hline\hline
\end{tabular}
}
\end{center}
\newpage
\begin{center}
{\emph{{\rm Table 6.} Numerical results for Experiment 2-II }\\
\scriptsize
\begin{tabular}{ l  l  l  l  l  l  l  l  l  l  l  l  l  l  l  l }
 \hline\hline
   Prob        &$n/m$           &Initial point          &NIO  &NII      &NF0       &NF        &FV      &CPU\\
\hline
 Svanberg-10   &10/30          &$(10,10,,\ldots,10)^T$   &3    &15      &19     &1320  &$15.731517$   &$0.28$\\
               &               &$(-10,-10,,\ldots,-10)^T$ &2   &16      &19     &1440  &$15.731517$  &$0.28$\\
\hline
 Svanberg-20   &20/60          &$(10,10,\ldots,10)^T$    &4    &22    &27     &4380    &$32.427932$   &$1.11$\\
               &               &$(-10,-10,,\ldots,-10)^T$ &3   &24     &28     &4560   &$32.427932$   &$1.23$\\
 \hline
Svanberg-30    &30/90          &$(10,10,,\ldots,10)^T$ &3 &25  &29        &6480       &$49.142526$    &$2.33$\\
               &               &$(-10,-10,,\ldots,-10)^T$&3  &24   &28    &6030       &$49.142526$    &$2.50$\\
\hline
Svanberg-40    &40/120        &$(10,10,,\ldots,10)^T$ &3    &28    &32    &9480      &$65.861140$    &$3.67$\\
               &              &$(-10,-10,,\ldots,-10)^T$ &3   &28    &32   &9480    &$65.861140$    &$3.58$\\
\hline
 Svanberg-50   &50/150       &$(10,10,\ldots,10)^T$  &14    &26    &41    &24150    &$82.581915$     &$7.84$\\
               &             &$(-10,-10,,\ldots,-10)^T$ &1  &34  &36    &12900      &$82.581912$     &$5.75$\\
 \hline
 Svanberg-80   &80/240       &$(10,10,\ldots,10)^T$  &2    &43    &46    &42720     &$132.749820$   &$17.09$\\
               &             &$(5,5,\ldots,5)^T$     &2    &47    &50    &63840     &$132.749824$   &$19.53$\\
 \hline
 Svanberg-100  &100/300     &$(10,10,\ldots,10)^T$  &3   &43    &47     &57300      &$166.197173$   &$26.81$\\
               &            &$(5,5,\ldots,5)^T$     &2   &62    &65    &112500      &$166.197178$   &$40.55$\\
 \hline
 Svanberg-150   &150/450    &$(10,10,\ldots,10)^T$  &40    &44   &85   &227700    &$249.818369$   &$130.41$\\
                &                  &$(5,5,\ldots,5)^T$ &3   &62    &66   &123750   &$249.818369$   &$96.16$\\
 \hline
 Svanberg-200   &200/600     &$(10,10,\ldots,10)^T$ &4    &78     &83   &219600    &$333.441310$   &$279.95$\\
                &            &$(5,5,\ldots,5)^T$    &2   &84   &87    &236400      &$333.441310$   &$287.45$\\
 \hline
 Svanberg-250   &250/750      &$(2,2,\ldots,2)^T$   &1    &85     &87    &275250    &$417.064989$   &$602.78$\\
                &             &$(3,3,\ldots,3)^T$   &1    &90     &92   &281250    &$417.064989$   &$593.48$\\
 \hline\hline
\end{tabular}
}
\end{center}

\textbf{Experiment 3.}  To show that ALG \ref{algorithm1} performs
cycle I for most of the iterations, we give Table 7 below relative
to Tables 1 and 2 for test problems. The row labeled $\sharp$ lists
the problem number as given in Tables 1 and 2, and the column
labeled $\sharp$ lists the number of cycle I and cycle II performed
by ALG \ref{algorithm1}, , i.e., N-cycle I and N-cycle II,
respectively.
 The results reported in Table 7 are encouraging.
Obviously, N-cycle I is much more than  N-cycle II.
 Especially, for problem 12, ALG \ref{algorithm1} always performs cycle I and
does not perform cycle II. This also illustrate Remark \ref{remark6}
from the viewpoint of numerical results. Thus, the cost of
computation for ALG \ref{algorithm1} is relatively small.
\newpage
\begin{center}
{\emph{{\rm {Table 7.}} Numerical results for performing Cycle I and Cycle II }\\
\scriptsize
\begin{tabular}{ c | c | c | c | c | c | c | c | c | c | c | c | c | c   }
 \hline\hline
 $\sharp$        &12  &29  &31  &33 &34  &35  &43  &44  &66  &76 &100  &113    &264      \\
\hline
  N-cycle I   &20  &8   &15  &9  &9   &6   &12  &7   &52    &17 &49   &11       &22            \\
 \hline
 N-cycle II  &0   &4   &2   &1  &6   &1   &2   &7   &12    &4  &8   &5        &1             \\
 \hline\hline
\end{tabular}
}
\end{center}
\section{Concluding remarks}\label{section6}
In this paper, we propose a new algorithm of combining (QP)
subproblem with SLE for solving nonlinear inequality constrained
optimization problems. The new algorithm starts from an arbitrarily
initial iteration point. In order to ensure the global convergence
of new algorithm, the search direction is obtained by a convex
combination of the master direction and an auxiliary direction,
which are solved by subproblem (QPs) and SLE (\ref{6}),
respectively. For overcoming the Maratos effect \cite{m1978}, a
higher-order direction is obtained by solving another SLE
(\ref{11}). Moreover, the iteration points can always enter into the
feasible set $X$ and only one SLE need to be solved after a finite
number of iterations. Using line search instead of arc search, our
new algorithm possesses global and superlinear convergence under
some mild assumptions without strict complementarity. Finally, some
numerical results show that new algorithm is promising.

As a further work of this paper, the techniques introduced in this
paper can be extended to solve general constrained optimization
problems and minimax problems.

%
%

\bibliographystyle{model1-num-names}
\bibliography{references_abbrname}
\end{document}